\documentclass{amsart}
\usepackage{amssymb,amsfonts}
\usepackage{ifthen}
\usepackage{epic,eepic}
\usepackage[all]{xy}

\newcommand\upl{\frac n2}
\newcommand{\sgn}{\operatorname{sgn}}

\newcommand\datver[1]{\def\datverp%
 {\par\boxed{\boxed{\text{Version: #1; Run: \today}}}}}
\datver{0.0a; Revised: 29-3-1997}

\newcommand\RR{\mathbb R}

\newcommand\pa{\partial}

\newcommand\dCI{\dot{\mathcal C}^{\infty}}

\newcommand\CIc{{\mathcal C}^{\infty}_{\text{c}}}

\newcommand\Lap{\varDelta}

\newcommand\Mfor{\text{ for }}

\newcommand\Min{\text{ in }}
\newcommand\Mon{\text{ on }}

\newcommand\restrictedto{\upharpoonright}
\newcommand{\WF}{\operatorname{WF}}

\newtheorem{theorem}{Theorem}

\theoremstyle{definition}

\theoremstyle{remark}

\newcommand\Cinf{\mathcal{C}^\infty}

\newcommand\Cminf{\mathcal{C}^{-\infty}}

\newcommand\E{\mathcal{E}}

\newcommand\loc{\operatorname{loc}}

\newcommand\clos{\operatorname{clos}}

\newcommand\const{\operatorname{const}}

\begin{document}

\title{Singularities and the wave equation on conic spaces}

\author{Richard B. Melrose} \author{Jared Wunsch} \address{Department of
Mathematics, MIT\\Cambridge MA 02139} \address{Department of Mathematics,
SUNY at Stony Brook\\Stony Brook NY 11794}


\maketitle

Introducing polar coordinates around a point in Euclidean space
reduces the Euclidian metric to the degenerate form
\begin{equation}\label{euclidian}
  dr^2 + r^2 \, d\omega^2
\end{equation}
where $r$ is the distance from the point and $d\omega^2$ is the round
metric on the sphere. If $X$ is an arbitrary manifold with boundary, the
class of \emph{conic metrics} on $X$ is modeled on this special
case. Namely, a conic metric is a Riemannian metric on the interior of $X$
such that for some choice of the defining function $x$ of the boundary ($x
\in \Cinf (X)$ with $\partial X = \{x=0\}$, $x \geq 0$, $d x \neq 0$ on
$\partial X$), the metric takes the form
\begin{displaymath}
  g=dx^2 + x^2 h \text{ on } X^\circ = X \backslash \partial X, 
\text{ near } \partial X. 
\end{displaymath}
Here $h$ is a smooth symmetric $2$-cotensor on $X$ such that
$h_0=h|_{\partial x}$ is a metric on $\pa X.$

In fact a general conic metric can be reduced to a form even closer to
\eqref{euclidian} in terms of an appropriately chosen product
decomposition of $X$ near $\partial X,$ that is, by choice of a smooth
diffeomorphism 
\begin{displaymath}
[ 0, \epsilon )_x \times \partial X \overset{F}{\longrightarrow} O
\subset X, \ O \text{ an open neighborhood of }\partial X \, .
\end{displaymath}
The normal variable in $x\in [0, \epsilon)$ is then a boundary defining
function, at least locally near $\partial X,$ and the slices $F
\restrictedto_{x=x_0}$ have given diffeomorphisms to $\pa X.$ Now such a
product decomposition can be chosen so that
\begin{equation}
  \label{RM}
  F^* g= dx^2 + x^2 h_x,\Min x<\epsilon ,
\end{equation}
where $h_x$ is a family of metrics on $\partial X$.

This reduced form is closely related to the behavior of geodesics near the
boundary. Up to orientation and parameterization there is a unique geodesic
reaching the boundary at a given point $p.$ In particular the normal
fibration of $X$ near $\partial X$ given by the segments $F([0,\epsilon )
\times\{p\})$, $p\in\partial X$, consists of geodesics which hit the
boundary, each at the corresponding point $p.$

We shall discuss here the behavior of solutions to the wave
equation
\begin{equation}
  \label{WE}
  (D^2_t - \Delta) u=0\text{ on } \RR \times X^\circ
\end{equation}
when $X$ is endowed with a conic metric, $\Delta$ is the associated
(positive) Laplacian on functions, and $D_t = -i \pa/\pa t$. For
simplicity we take $X$ to be compact. It is only really important that
$\partial X$ be compact.

Our primary concern is to describe the phenomenon of the propagation of
singularities for solutions to \eqref{WE}. To do so it is necessary to
understand the behavior of solutions in a way related to the functional
analytic domain of $\Delta.$ For the moment we simply say that we are
dealing with `admissible' solutions. This condition is explained further
below.

In the interior of $X$ the propagation of singularities, described
precisely in terms of the notion of wavefront set, was treated in detail
by H\"ormander (\cite{Hormander9}). We paraphrase H\"ormander's result here as
\begin{quotation}
\label{HT}
``Singularities travel along null bicharacteristics, which in the case of
the wave equation project to time-para\-meter\-ized geodesics.''
\end{quotation}
Thus, in the microlocal sense of singularities described by the
wave front set, a \emph{bicharacteristic segment}, which covers
a light ray, either consists completely of singularities for a given
solution or the solution has no singularity along it.

This quite adequately describes the propagation of singularities except
where a light ray hits the boundary at some point and at some time. Here a
`splitting' of singularities will usually take place. This is generally
called a \emph{diffractive} effect. The contrapositive of this effect can
be succinctly stated as follows:
\begin{quotation}
``If no singularity reaches the boundary at time $\bar t$ then no singularity
leaves at time $\bar t.$''
\end{quotation}

The point here is that the regularity along any one of the `radial' rays
leaving the boundary at a given time is related, in general, to the
singularities on all the incoming rays (although there are two separate
components, as described below) arriving at the boundary at that
time. Thus, even if singularities arrive at the boundary at time $\bar t$ along
just one ray, they will in general depart along all rays leaving the
boundary at time $\bar t.$

There are, however, some important exceptions to this general spreading of
singularities. For instance let $X$ be a conic manifold with `trivial' conic
metric defined by the blowup of a point in a smooth Riemannian manifold. In
this case, because of H\"ormander's theorem on interior singularities, the
singularities are carried outward only on the one ray continuing the
incoming ray in the original manifold.

For a general conic metric there is a similar notion of the `geometric
continuation' of an incoming geodesic which hits the boundary. For a
trivial conic metric obtained from a blowup, the boundary metric $h_0$ is the
standard metric on the sphere. The geometrically related incoming and
outgoing rays hit this sphere at antipodal points; these can also be
thought of as the points separated by geodesics of length $\pi$ on the unit
sphere. In the case of a general conic metric we mimic this by defining the
relation
\begin{multline}
  \label{GS} G(p) = \\
\{ q \in \partial X ;\,\, \exists \text{ a geodesic in
  } \partial X \Mfor h_0 \text{ of length } \pi \text{ with end points }p,q\}.
\end{multline}
In general of course, $G(p)$ is not smooth, but generically it is a
hypersurface with Lagrangian singularities; it is always the projection of
a smooth Lagrangian relation.

A geometric refinement of the diffraction result is obtained by considering
the order of singularity with respect to Sobolev spaces and an additional
`second microlocal' regularity condition. For simplicity suppose that the
(admissible) solution $u$ is singular only near $\partial X$ and only near
a single incoming ray hitting the boundary at time $\bar{t}$ and at the
point $p.$ In the past (for $t<\bar t$) we may suppose that the solution is
locally in some Sobolev space $H^s.$ Suppose further that the singularities
of the solution are not too strongly focused on $\pa X$ insofar as
\emph{tangential smoothing} raises the overall regularity, that is, for
some $k, \ell >0$,
\begin{equation}
  \label{TR}
  (\Delta_0 +1)^{-k} u \in H^{s+\ell}_{\loc}
  \text{ in } t < \bar{t} \text{ near }\pa X.
\end{equation}
Under these two assumptions and the additional requirement that
\begin{equation}
  \label{ER}
  0 < \ell <\upl,
\end{equation}
we obtain the following `geometric theorem'.

\begin{quotation}
``If an admissible solution is singular only near an incoming ray arriving
at $\partial X$ at time $\bar{t}$ and \eqref{TR} and \eqref{ER} hold, then
on outgoing rays with initial point in the complement of $G(p)$,
\begin{equation}
  \label{OR}
u \in H^{s+\ell-\frac 12 -\epsilon}_{\loc} \ \forall\ \epsilon >0
   \text{ in } t > \bar{t} \text{ near } \partial X.
\end{equation}
\end{quotation}
\begin{figure}
\setlength{\unitlength}{0.0003in}
\begingroup\makeatletter\ifx\SetFigFont\undefined%
\gdef\SetFigFont#1#2#3#4#5{%
  \reset@font\fontsize{#1}{#2pt}%
  \fontfamily{#3}\fontseries{#4}\fontshape{#5}%
  \selectfont}%
\fi\endgroup%
{\renewcommand{\dashlinestretch}{30}
\begin{picture}(7520,5964)(0,-10)
\put(1133,3012){\ellipse{2250}{3900}}
\blacken\path(2266.230,4144.484)(2033.000,4062.000)(2277.607,4025.024)(2266.230,4144.484)
\path(2033,4062)(3608,4212)(6008,4437)
\path(2183,2337)(4883,2037)(6758,1812)
\blacken\path(6512.561,1781.022)(6758.000,1812.000)(6526.858,1900.167)(6512.561,1781.022)
\path(1958,1662)(4283,1287)(7058,837)
\blacken\path(6811.490,816.191)(7058.000,837.000)(6830.699,934.643)(6811.490,816.191)
\path(2258,2937)(4658,2712)(6233,2562)
\blacken\path(5988.393,2525.024)(6233.000,2562.000)(5999.770,2644.484)(5988.393,2525.024)
\path(1133,4962)(1134,4962)(1137,4962)
	(1142,4962)(1150,4963)(1163,4964)
	(1179,4964)(1200,4966)(1226,4967)
	(1257,4969)(1294,4971)(1337,4973)
	(1385,4976)(1438,4979)(1497,4982)
	(1560,4986)(1628,4990)(1700,4994)
	(1775,4999)(1854,5004)(1935,5008)
	(2018,5014)(2102,5019)(2187,5024)
	(2273,5030)(2359,5035)(2445,5041)
	(2530,5047)(2614,5053)(2697,5058)
	(2778,5064)(2858,5070)(2936,5076)
	(3012,5082)(3086,5087)(3158,5093)
	(3229,5099)(3298,5105)(3364,5111)
	(3429,5117)(3493,5123)(3554,5129)
	(3614,5135)(3673,5141)(3730,5147)
	(3787,5154)(3842,5160)(3896,5166)
	(3949,5173)(4002,5180)(4054,5187)
	(4105,5194)(4156,5201)(4207,5209)
	(4258,5217)(4308,5224)(4360,5233)
	(4412,5241)(4464,5250)(4517,5259)
	(4569,5269)(4622,5279)(4676,5289)
	(4729,5299)(4784,5310)(4840,5321)
	(4896,5333)(4953,5345)(5012,5358)
	(5072,5371)(5134,5385)(5197,5400)
	(5261,5415)(5327,5430)(5395,5446)
	(5465,5463)(5537,5481)(5610,5499)
	(5685,5517)(5762,5536)(5840,5556)
	(5919,5576)(6000,5597)(6081,5618)
	(6163,5639)(6245,5660)(6327,5682)
	(6409,5703)(6489,5724)(6568,5745)
	(6644,5765)(6718,5785)(6789,5804)
	(6856,5822)(6918,5838)(6977,5854)
	(7030,5868)(7078,5881)(7120,5893)
	(7157,5903)(7189,5911)(7215,5919)
	(7237,5924)(7253,5929)(7265,5932)
	(7274,5934)(7279,5936)(7282,5937)(7283,5937)
\path(1133,1062)(1134,1062)(1137,1062)
	(1143,1061)(1151,1060)(1164,1059)
	(1181,1057)(1203,1055)(1230,1053)
	(1263,1050)(1302,1046)(1346,1042)
	(1396,1037)(1452,1032)(1513,1026)
	(1580,1019)(1651,1012)(1726,1005)
	(1805,997)(1887,989)(1971,981)
	(2058,972)(2146,964)(2235,955)
	(2325,946)(2415,937)(2504,927)
	(2593,918)(2681,909)(2767,900)
	(2852,891)(2935,882)(3017,874)
	(3096,865)(3174,856)(3249,848)
	(3323,840)(3394,831)(3464,823)
	(3532,815)(3598,807)(3662,799)
	(3725,791)(3786,784)(3846,776)
	(3904,768)(3961,760)(4018,752)
	(4073,744)(4128,736)(4182,728)
	(4235,720)(4288,712)(4341,704)
	(4393,696)(4446,687)(4500,678)
	(4554,669)(4608,660)(4662,650)
	(4716,640)(4771,630)(4826,620)
	(4882,610)(4938,599)(4995,588)
	(5054,576)(5113,564)(5173,552)
	(5235,539)(5299,525)(5364,512)
	(5430,497)(5498,482)(5568,467)
	(5640,451)(5714,435)(5789,418)
	(5866,400)(5945,382)(6025,364)
	(6107,345)(6190,326)(6274,306)
	(6358,287)(6442,267)(6527,247)
	(6610,227)(6693,208)(6773,189)
	(6852,170)(6928,152)(7001,134)
	(7069,118)(7134,103)(7193,88)
	(7248,75)(7297,63)(7341,53)
	(7379,43)(7411,35)(7439,29)
	(7460,24)(7477,19)(7490,16)
	(7498,14)(7504,13)(7507,12)(7508,12)
\put(1783,4137){$p$}
\put(3908,3987){\makebox(0,0)[lb]{\smash{{{\SetFigFont{10}{14.4}{\rmdefault}{\mddefault}{\updefault}incoming singularity}}}}}
\put(3683,2112){\makebox(0,0)[lb]{\smash{{{\SetFigFont{10}{14.4}{\rmdefault}{\mddefault}{\updefault}outgoing regularity}}}}}
\put(1158,2412){${G(p)}^\complement$}
\end{picture}
}
\caption{Regularity occurs at outgoing rays from points not $\pi$-related to $p$.}
\end{figure}
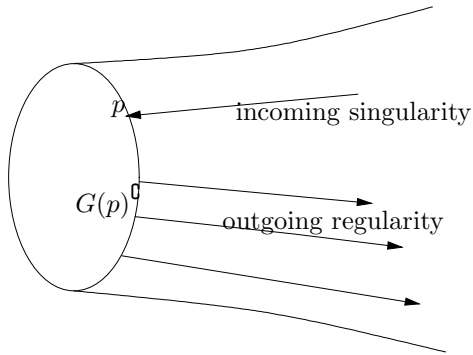

When slightly generalized, as described below, this result applies to the
fundamental solution $$\frac{\sin t{\sqrt\Lap}}{\sqrt\Lap}$$ with pole
close to $\partial X$ and with $\ell < \frac{n-1}{2}$. The diffractive
theorem merely tells us that if the pole is specified at $(\bar x,p)$ at
$t=0$, then singularities cannot emanate from $\partial X$ except at time
$t=\bar x$.  On the other hand, while `strong' singularities can emanate
from all points in $G(p),$ the geometric theorem tells us that the solution
is microlocally more regular on rays starting from $\pa X$ at $t=\bar x$
but with initial point outside $G(p).$

In the special case in which the metric $g$ takes precisely the `product' form
\begin{equation}\label{productmetric}
g = dx^2 + x^2 h(y, dy),
\end{equation}
near the boundary, Cheeger and Taylor \cite{Cheeger-Taylor2,
Cheeger-Taylor1} have given an explicit analysis of the fundamental
solution constructed by separation of variables.  (See also the discussion
by Kalka-Menikoff \cite{Kalka-Menikoff1}.) They show a stronger form of the
result above, including a more precise regularity estimate.  A version of
the results of Cheeger-Taylor has been established in the analytic category
by Rouleux \cite{Rouleux}. Lebeau \cite{Lebeau4, Lebeau5} has also
obtained a diffractive theorem in the setting of manifolds with corners in
the analytic category.

The first author acknowledges partial support from NSF grant DMS9625714.
The second author was partially supported by an NSF VIGRE instructorship at
Columbia University and by an NSF postdoctoral fellowship, and is grateful
to Alberto Parmeggiani and Cesare Parenti for helpful conversations.

\section{Friedrichs extension}
\label{sec:1}

To describe the admissibility condition, near the boundary, for
solutions to \eqref{WE} we first describe the domain of the
Laplacian for a conic metric. We take the Friedrichs extension
of $\Lap$. By definition, $\Lap$ is associated to the Dirichlet form
\begin{equation}
  \label{DF}
  F(u,v) = \int_X \langle du,dv \rangle_g \, dg \,,\,\,
  u,v \in \Cinf_c (X^\circ)
\end{equation}
Hence, $dg$ is the metric volume form; in this case
\begin{displaymath}
  dg = \varphi x^{n-1} \, dx \, dh_0 \text{ near }
  \partial X, \quad n=\dim X \, , \, \varphi \in \Cinf ,
  \varphi >0 \, .
\end{displaymath}
The inner product in \eqref{DF} is that induced, by duality, by
the metric on $T^* X^\circ.$ Following Friedrichs we define, 
\begin{equation*}
D (\Lap^{1/2}) = \clos\left\{\CIc(X^\circ) \text{ w.r.t. }
F(u,u)+\|u\|^2_{L^2_g} \right\},
\end{equation*}
whenever $X$ is a compact conic manifold with boundary of dimension $n \geq
2.$

This is a Hilbert space with dense injection $D (\Lap^{1/2})
\hookrightarrow L^2_g (X)$ so there is a dual injection $L^2_g(x)
\hookrightarrow (D(\Lap^{1/2}))'.$ The natural operator $\Lap:D(\Lap^{1/2})
\to (D(\Lap^{1/2}))'$ is determined  by
\begin{displaymath}
  (\Lap u,\varphi)_{L^2_g} = F (u,\varphi) \
  \forall\ u,\varphi \in D(\Lap^{1/2}).
\end{displaymath}
Then the Friedrichs extension of $\Lap$ is the unbounded operator
with domain
\begin{displaymath}
  D(\Lap) = \left\{u\in D(\Lap^{1/2}), \,\, \Lap u \in L^2_g (x)\right\} \, .
\end{displaymath}
In this case, it is a self-adjoint, non-negative operator with discrete
spectrum of finite multiplicity. This allows its complex powers to be
defined by reference to an eigenbasis. The real powers are isomorphisms off
the null space, which consists precisely of the constants. Each of the
powers is therefore a Fredholm map
\begin{displaymath}
  \Delta^s : D(\Lap^s) \to L^2_g (X)\ \forall\  s \in \RR.
\end{displaymath}
with null space the constants and range the orthocomplement of the constants.
The domains form a scale of Hilbert spaces, and
\begin{displaymath}
  D(\Lap^s) \hookrightarrow D (\Lap^t) \text{ is dense }
  \forall\ s \geq t
\end{displaymath}
with $D(\Lap^0)=L^2_g(X)$.

For our purposes it is also important to note when the domains consist of
extendible distributions, i.e.\ those dual to $\dCI(X)$.  This is the case
only for $s> -\frac{n}{4}$ and more precisely
\begin{equation}
  \label{DI}
  \dCI (X) \hookrightarrow D(\Lap^s) \hookrightarrow \Cminf (X)
\end{equation}
are dense inclusions for $-\frac{n}{4} < s < \frac{n}{4}.$ The limits of
this range correspond to the occurrence of formal solutions of $\Lap u=0.$

\section{Wave group}
\label{sec:2}

The Cauchy problem for the wave equation
\begin{equation}\label{wave-equation}
\begin{gathered}
  (D^2_t -\Delta) u = 0,\Mon \RR \times X^\circ \\
  u|_{t=0} = u_0, \ D_t u|_{t=0} =u_1
\end{gathered}
\end{equation}
has a unique solution
\begin{eqnarray*}
  u \in C^0(\RR ; D (\Lap^{1/2})) \cap   C^1 (\RR ; L^2_g (x))\\
  \forall\ (u_0,u_1) \in \E = \E_1 = D(\Lap^{1/2}) \oplus L^2_g (X) \, .
\end{eqnarray*}
These `finite energy solutions' are the main object of study
here. More generally, with the equation interpreted in $\Cminf (\RR ; 
D (\Lap^{\frac{s}{2}-1}))$ the Cauchy problem has a unique solution
\begin{equation}
  \label{CPS}
\begin{gathered}
  u \in C^0(\RR ; D (\Lap^{\frac s2}))\cap C^1(\RR;D(\Lap^{\frac s2-\frac12}))\\
\forall\ (u_0,u_1) \in \E_s=
D (\Lap^{\frac s2}) \oplus D (\Lap^{\frac s2-\frac12}) \, .
\end{gathered}
\end{equation}
The regularity hypothesis on the solution can be weakened to
\begin{equation*}
u \in L^2_{\loc} (\RR ; D(\Lap^{\frac{s}{2}})) \cap
  H^1_{\loc} (\RR ; D (\Lap^{\frac{s}{2}-\frac{1}{2}}))
\end{equation*}
without changing the unique solvability.

Notice that these calculations are consistent under decrease of $s.$
Furthermore, partial hypoellipticity in $t$ shows that the solution to
\eqref{wave-equation} satisfies
\begin{equation}
  \label{PH}
  u \in H^{-k}_{\loc} (\RR ; D (\Lap^{\frac{s}{2}+\frac{k}{2}}))
  \quad \forall\ k \in \RR \, .
\end{equation}
An admissible solution to the wave equation is one that satisfies
\begin{equation}
  \label{AS}
  u \in H^{p}_{\loc} (\RR ; D(\Lap^{\frac{q}{2}}))
  \text{ for some } q,p \in \RR
\end{equation}
with \eqref{WE} holding in $H^{-p-2}_{\loc} (\RR ; D (\Lap^{\frac{q}{2}-1}))$.
Such a solution automatically satisfies \eqref{PH} for some $s.$

These statements can be reinterpreted in terms of the wave group 
\begin{equation}
  \label{WG}
  U(t) :\binom{u_0}{u_1} \mapsto \binom{u(t)}{D_t u(t)},\
  U (t) : \E_s \to \E_s \ \forall\ s.
\end{equation}

\section{H\"ormander's theorem}
\label{sec:3}

Let $M$ be a manifold without boundary. The wave front set of a
distribution $u \in \Cminf (M)$ is a closed subset of the cosphere bundle
\begin{displaymath}
  \WF (u) \subset S^*M.
\end{displaymath}
It may be defined by decay properties of the localized Fourier transform,
or the FBI (Fourier--Bros--Iagonitzer) transform, or by testing with
pseudodifferential operators. The projection $\pi (\WF(u)) \subset M$ is
exactly the $\Cinf$ singular support, the complement of the largest open subset
of $M$ to which $u $ restricts to be $\Cinf$.

A refined notion of wavefront set is the Sobolev-based wavefront set,
denoted $\WF_s$; this is a closed subset of $S^* M$, where now the
projection is the complement of the largest open subset of $M$ to which $u$
restricts to be $H^s$.

If $u$ satisfies a linear differential equation, $Pu=0,$ then
\begin{displaymath}
  \WF (u) \subset \Sigma (P) \subset S^* M 
\end{displaymath}
when $\Sigma (P)$ is the characteristic variety of $P$, the set
on which its (homogeneous) principal symbol, $p,$ vanishes.

If $p$ is real then the symplectic structure on $T^* M$, or the contact
structure on $S^* M,$ defines a `bicharacteristic' direction field $V_P$ on
$S^*M$, tangent to $\Sigma (P).$ The integral curves of $V_P$ are called
bicharacteristics; those lying in $\Sigma (P)$ are called null
bicharacteristics.

\begin{theorem}[H\"ormander]
  \label{th:Hormander} Let $P$ be a (pseudo)-differential operator with
  real principal symbol. If $Pu=0$ then $WF (u) \subset \Sigma (P)$ is a
  union of maximally extended null bicharacteristics.
\end{theorem}
The same result also holds with $\WF$ replace by $\WF_s$ for any $s$.

In our case, $M=\RR \times X^\circ$ so $T^* M = T^* \RR \times T^*
X^\circ$. The principal symbol of the d'Alembertian is $\tau^2 - | \cdot
|^2_g$, where $\tau$ is the dual variable to $t$ and $| \cdot
|_g$ is the (dual) metric on $T^* X^\circ$. Then 
\begin{displaymath}
  \Sigma (P) = \Sigma_+ (P) \cup \Sigma_- (P) \subset S^*M
\end{displaymath}
when $\Sigma_I (P) \cong \RR \times S^* X^\circ$ are the disjoint parts of
$\Sigma (P)$ in $\tau >0$ and $\tau <0$. In this representation of $\Sigma
(P)$ the null bicharacteristics are geodesics on $X^\circ$, lifted
canonically to $S^*X^\circ$, with $t$ as affine parameter. Thus, for the
wave equation over $X^\circ$, H\"ormander's theorem does indeed reduce to
the informal propagation statement desribed above.

Combined with standard results relating the singularities of the solution
to singularities of the initial data, H\"ormander's theorem applied to
the wave equation on a conic manifold to yields complete information on the
behavior of singularities except along bicharacteristics lying above
geodesics which hit the boundary.

\section{Diffractive theorem}
\label{sec:4}
On parametrized geodesic segments with an end point on the boundary, the
defining function $x$ is either strictly increasing or strictly decreasing
near the boundary. For each sign of $\tau$ and for each $\bar t \in \RR$
the bicharacteristics covering such geodesics which hit the boundary at
$t=\bar t$ and along which $t$ is increasing (resp.\ decreasing) as $x$
decreases, form a smooth submanifold of $\{t< \bar t\} \times X$ (resp.\ $\{t>
\bar t\} \times X).$  We denote these `radial' surfaces (near $\pa X$) by
\begin{displaymath}
  R_{\pm,I} (\bar t) \text{ and } 
  R_{\pm,O} (\bar t) \subset \Sigma (P)
\end{displaymath}
where $\pm$ is the sign of $\tau$ and $I,O$ refers to whether
these are `incoming' or `outgoing' and hence, equivalently, whether 
they lie in $t<\bar t$ or $t>\bar t.$

\begin{theorem}[Diffractive theorem]
  \label{th:Diffraction}
  If $u$ is an admissible solution to \eqref{WE} then for any $\bar t 
  \in \RR$, $s\in \RR$, $\sigma=\pm$,
  \begin{equation*}
R_{\sigma ,I} (\bar t) \cap \WF_s (u) = \emptyset
    \Rightarrow R_{\sigma ,O} (\bar t) \cap \WF_s (u) =\emptyset.
  \end{equation*}
\end{theorem}

Here, $\WF_s (u)$ is the wave front set computed relative to the
Sobolev space $H^s,$ locally in the interior.

This is a precise form of the diffractive result described informally
above. Notice that the singularities for different signs of $\tau$ are
completely decoupled. This does not, however, represent any refinement in
terms of propagation along the underlying geometric rays, since all
geodesics are covered by bicharacteristics with $\tau$ fixed and of either
sign.

The proof of this result is discussed briefly below in \S\ref{sec:7}.

\section{Geometric theorem}
\label{sec:5}

\begin{figure}
\input{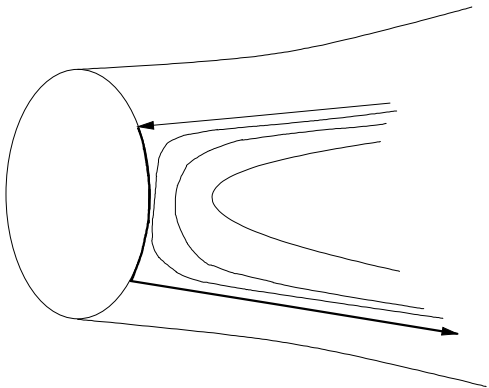}
\caption{A sequence of geodesics nearly missing the boundary, and the three
segments to which they limit.\label{fig2}}
\end{figure}

Consider a geodesic on $X$ which hits the boundary at a point $p\in\pa X.$
An open set of perturbations of the geodesic, meaning geodesics starting
near some interior point on the geodesic and with initial tangent close to
the tangent to the geodesic, will miss the boundary. A limit of such curves
as the perturbation vanishes consists of three segments. The first is the
incoming geodesic segment. The second is a geodesic segment in the
boundary, of length $\pi$. The third is the outgoing geodesic from the end
point of the boundary segment, which is therefore a point in $G(p)$ as
defined in \eqref{GS} (see Figure~\ref{fig2}).
Thus it is reasonable to suppose that,
amongst the outgoing bicharacteristics leaving the boundary at time $\bar
t,$ those with initial points in $G(p)$ will be more closely related to an
incoming bicharacteristic with end point $p$ arriving at time $\bar t.$ We
call these the geometrically-related bicharacteristics (or geodesics).

For instance, if there are incoming singularities on a single ray the
singularities on the `non-geometrically-related' outgoing bicharacteristics
might be expected to be weaker than the incoming singularity. However, this
is \emph{not} in general the case. To obtain such a geometric refinement of
the diffraction result we need to impose an extra `nonfocusing' assumption.

\begin{theorem}[Geometric theorem]
Let $u $ be an admissible solution to \eqref{WE} and let $\sigma=\pm$.
Suppose that $R_{\sigma,I}(\bar t) \cap \WF_s u = \emptyset$ near $\pa X$.
Suppose additionally that for some $k$
and $0<\ell <\upl$
\begin{equation}
  \label{TM}
  \WF_{s+\ell}(1+\Lap_{\pa X})^{-k} \, u \cap  R_{\sigma,I}(\bar t) = \emptyset.
\end{equation}
For any $0<r<\ell-1/2,$ if no incoming bicharacteristic hitting the
boundary at time $\bar t$ at a point in $G(p)$ with $\sgn \tau = \sigma$ is
in $\WF_{s+r} u$, then the outgoing bicharacteristic with initial point $p
\in \pa X$ and $\sgn \tau=\sigma$ is not in $\WF_{s+r} u$ either.
\end{theorem}

Thus the additional assumption \eqref{TM} allows regularity on
outgoing rays to be deduced from regularity in the incoming
geometrically-related rays up to the corresponding level above
`background' regularity.

As already noted, this result may be applied to the fundamental solution
with initial point near the boundary. If the initial pole of the
fundamental solution is sufficiently close to the boundary then there is a
unique short geodesic segment from it to the boundary, arriving at a point
$p.$ If $\bar t$ is the length of the segment then, provided $\bar t$ is
small enough, \eqref{TM} holds with $s<-\frac{n}{2}+1$ for any $\ell <
\frac{n-1}{2}$. It follows that on $R_{\pm ,O}(\bar t),$ the outgoing set,
the fundamental solution is in $H^{-\epsilon}$, for all $\epsilon >0$,
microlocally near the non-geometrically related rays, those with end point
not in $G(p)$, whereas the general regularity is
$H^{-\frac{n}{2}+1-\epsilon}$ for all $\epsilon >0$. This is a gain of
`nearly' $\frac{n}{2}-1$ derivatives over the background regularity.

In this way we extend part of the result of Cheeger and Taylor
\cite{Cheeger-Taylor2,Cheeger-Taylor1} in the product case
\eqref{productmetric} to the general conic case.  Note, however, that
inspection of the fundamental solution constructed in
\cite{Cheeger-Taylor2} reveals the diffracted wave to be $(n-1)/2-\epsilon$
derivatives smoother than the incident wave; here, we only obtain
$(n-2)/2-\epsilon$ derivatives of improvement, hence our results are not
sharp in the product case.  Even in the non-product case, we conjecture
that the better estimate should hold.

\section{Spherical conormal waves}
\label{sec:6}\

Around a given point $q$ in a compact Riemann manifold there are
`spherical' conormal waves which are singular only on the spherical
surfaces $r=\pm t$, for small $t$ of both signs. These just correspond to
conormal data at $t=0$ at the (fictive) cone point $q$. An important
example is the fundamental solution, in which case the result follows from
Hadamard's construction. In the more general case of a conic manifold with
boundary there are similar contracting, and then expanding, conormal waves.

\begin{theorem}
  If $u$ is an admissible solution near $\partial X$ and $t=0$
  which is conormal to $t=-x$ for $t<0$ then it is conormal to
  $t=x$, near the boundary, for small $t>0$.
\end{theorem}

These conormal solutions to the wave equation in the general conic case are
at the opposite extreme to those considered in the Geometric Theorem
above. Namely, they are already smooth in the tangential variables, so no
tangential smoothing in the sense of \eqref{TR} is possible. Further analysis of the
structure of these waves shows that the principal symbols undergo a
transition at $x=0$, the boundary, given by the scattering matrix for the
model cone with the same boundary metric.  Since this scattering matrix
should have full support in general, this provides counterexamples to any
extension of the geometric theorem in which the tangential smoothing
condition is dropped.

\section{Methods}
\label{sec:7}

The basic method we use is microlocal, but non-constructive. It is a direct
extension of one of the proofs by H\"ormander of the interior propagation
theorem. This `positive' commutator method is itself a microlocalization of
the energy method for hyperbolic equations.  In it a `test'
pseudodifferential operator, $A$, is applied to the equation and the
essential positivity of the symbol of the commutator $\frac{1}{i}[P,A]$
gives a local regularity estimate on the solution.

To extend this method to cover behavior of solutions near the boundary we
replace the ordinary notions of wavefront set, pseudodifferential operators
and microlocalization with versions appropriately adapted to the
geometry. When considering the Laplacian itself on the manifold with
boundary with conic metric, the appropriate notion is that of a weighted
$b$-pseudodifferential operator (see \cite{MR96g:58180}). This for instance
allows the precise description of the domains of the powers of $\Lap$ which
is used at various points in the argument.

However, for the wave operators for the conic Laplacian the appropriate
notion corresponds to the `edge' calculus of pseudodifferential operators
discussed originally by Mazzeo \cite{Mazzeo4}, arising from a filtration of
the boundary (see also Schulze \cite{Schulze2}). In this case, the manifold
with boundary is $X \times \RR$ and the fibers of the boundary are the
surfaces $t=\const$. Thus $t$ is the base variable of the fibration.

To the edge calculus of pseudodifferential operators, given by
microlocalization from the differential operators generated by $xD_x$,
$D_y$ (where the $y$'s are tangential variables) and $xD_t$, we associate a
notion of wavefront set.  We can prove the propagation
theorem analogous to that of H\"ormander for this `edge' wavefront set. However, in this
new sense, $D^2_t-\Lap$ is \emph{not} globally of principal type but rather
has two radial surfaces. These correspond to the end points of
bicharacteristics arriving at, and leaving from, the boundary. At these
surfaces there are restrictions on the propagation results, very closely
related to those for scattering Laplacians in \cite{Melrose43}.

These propagation estimates form the basis of both the diffractive and
geometric theorems. In the former we combine the estimates with a variant
of the one-dimensional FBI transform, scaled with respect to the normal
variable $x$. This reduces the diffractive result to an iterative
application of a uniqueness theorem for the Laplacian on the model,
non-compact cone.

To obtain the geometric theorem, showing that the outgoing singularities on
non-geometrically related rays are weaker than the incoming ones, we use a
division theorem. The additional hypothesis of microlocal tangential
smoothing is shown to imply that the solution actually lies in a weighted
Sobolev space with a higher $x$ weight (hence more `divisible' by $x$) than
is given, \emph{a priori}, by energy conservation. This allows the
microlocal propagation results indicated above to be pushed further at the
outgoing radial surface and so yields the extra regularity.

\section{Applications and extension}
\label{sec:8}

The propagation of singularities results of the type discussed above should
allow estimates of the spectral counting function as shown originally by
Ivrii (\cite{Ivrii2}, see also \cite{MR86f:35144} and \cite{Hormander3}).

We expect these methods to extend to more complicated geometries, 
including manifolds with corners and iterated conic spaces.

\bibliography{all,My}

\providecommand{\bysame}{\leavevmode\hbox to3em{\hrulefill}\thinspace}
\begin{thebibliography}{10}

\bibitem{Cheeger-Taylor2}
Jeff Cheeger and Michael Taylor, \emph{On the diffraction of waves by conical
  singularities. {I}}, Comm. Pure Appl. Math. \textbf{35} (1982), no.~3,
  275--331, MR84h:35091a.

\bibitem{Cheeger-Taylor1}
\bysame, \emph{On the diffraction of waves by conical singularities. {I}{I}},
  Comm. Pure Appl. Math. \textbf{35} (1982), no.~4, 487--529, MR84h:35091b.

\bibitem{Hormander3}
L.~H{\"o}rmander, \emph{The analysis of linear partial differential operators},
  vol.~3, Springer-Verlag, Berlin{,} Heidelberg{,} New York{,} Tokyo, 1985.

\bibitem{Hormander9}
Lars H{\"o}rmander, \emph{On the existence and the regularity of solutions of
  linear pseudo-differential equations}, Enseignement Math. (2) \textbf{17}
  (1971), 99--163.

\bibitem{Ivrii2}
V.~Ivrii, \emph{On the second term in the spectral asymptotics for the
  {L}aplace-{B}eltrami operator on a manifold with boundary}, Funct. Anal.
  Appl. \textbf{14} (1980), 98--106.

\bibitem{Kalka-Menikoff1}
M.~Kalka and A.~Menikoff, \emph{The wave equation on a cone}, Comm. Partial
  Differential Equations \textbf{7} (1982), no.~3, 223--278, MR83j:58110.

\bibitem{Lebeau4}
G.~Lebeau, \emph{Propagation des ondes dans les vari\'et\'es \`a coins},
  S\'eminaire sur les \'Equations aux D\'eriv\'ees Partielles, 1995--1996,
  \'Ecole Polytech., Palaiseau, 1996, MR98m:58137, pp.~Exp. No. XVI, 20.

\bibitem{Lebeau5}
Gilles Lebeau, \emph{Propagation des ondes dans les vari\'et\'es \`a coins},
  Ann. Sci. \'Ecole Norm. Sup. (4) \textbf{30} (1997), no.~4, 429--497,
  MR98d:58183.

\bibitem{Mazzeo4}
R.~Mazzeo, \emph{Elliptic theory of differential edge operators {I}}, Comm. in
  P.D.E. \textbf{16} (1991), 1615--1664.

\bibitem{Melrose43}
R.B. Melrose, \emph{Spectral and scattering theory for the {L}aplacian on
  asymptotically {E}uclidian spaces}, Spectral and scattering theory (Sanda,
  1992) (M.~Ikawa, ed.), Marcel Dekker, 1994, pp.~85--130.

\bibitem{MR86f:35144}
Richard Melrose, \emph{The trace of the wave group}, Microlocal analysis
  (Boulder, Colo., 1983), Amer. Math. Soc., Providence, R.I., 1984,
  pp.~127--167.

\bibitem{MR96g:58180}
Richard~B. Melrose, \emph{The {A}tiyah-{P}atodi-{S}inger index theorem}, A K
  Peters Ltd., Wellesley, MA, 1993.

\bibitem{Rouleux}
Michel Rouleux, \emph{Diffraction analytique sur une vari\'et\'e \`a
  singularit\'e conique}, Comm. Partial Differential Equations \textbf{11}
  (1986), no.~9, 947--988.

\bibitem{Schulze2}
B.-W. Schulze, \emph{Boundary value problems and edge pseudo-differential
  operators}, Microlocal analysis and spectral theory (Lucca, 1996), Kluwer
  Acad. Publ., Dordrecht, 1997, pp.~165--226.

\end{thebibliography}
\bibliographystyle{amsplain}

\end{document}